\documentclass[12pt]{amsart}
\usepackage{amsmath,amsfonts,amsthm,amscd,graphicx,amssymb,verbatim}
\usepackage{enumerate}
\usepackage{tikz}
\usetikzlibrary{arrows}
\usepackage{euscript}
\usepackage{soul}

\usepackage{color}

\newcommand{\op}{\mathit{op}}
\theoremstyle{plain}
\newtheorem{thm}{Theorem}[section]

\newtheorem{lem}[thm]{Lemma}

\theoremstyle{definition}

\newtheorem{ex}[thm]{Example}

\newtheorem{defn}[thm]{Definition}

\newcommand{\maps}{\ensuremath{\colon}}

\newcommand{\obj}{\operatorname{obj}}

\newcommand{\skc}{\operatorname{sk}^{c}}
\DeclareMathOperator{\colim}{colim}

\begin{document}

\title{\bf A Twisted Version of the Classifying Space Functor}
\author{Asl\i \ G\"{u}\c{c}l\"{u}kan \.{I}lhan and \"{O}zg\"{u}n \"{U}nl\"{u}}

\address{Asli G\" u\c cl\" ukan \.{I}lhan\\
Department of Mathematics\\
 	Dokuz Eyl\"ul University\\  \.Izmir, Turkey}

\email{asli.ilhan@deu.edu.tr}

\address{\"{O}zg\"{u}n \"{U}nl\"{u} \\
Department of Mathematics \\
Bilkent University \\ Ankara, Turkey}

\email{unluo@fen.bilkent.edu.tr}

\thanks{The authors are supported by T\"UB\.ITAK-TBAG/117F085.}

\begin{abstract}
It is known that there is a weak-equivalence between the geometric realization of a simplicially enriched small category  and its cofibrant replacement \cite{Riehl}. In this paper, we show that when only small categories are considered there exists a homeomorphism between these geometric realizations. We also discuss the naturality of these homoemorphisms. The inclusion of the category of small categories to the category of simplicially enriched categories, the cofibrant replacement of simplicially enriched categories, and the geometric realization of simplicially enriched categories are three composable functors. Hence one can ask if the collection of all these homeomorphisms gives a natural transformation from the composition of these three functors to the classifying space functor. We show that this is almost the case and that this composition can be considered as some twisted version of the classifying space functor.
\end{abstract}

\maketitle
\section{Introduction}
\label{sect:introduction}

One way of constructing a topological space with a group action is to use a gluing data which is given up to homotopy. A gluing data may be given by a homotopy diagram which is a simplicially enriched functor from a cofibrant replacement of a simplicially enriched category to $\mathbf{Top}$, the category of topolgical spaces. In \cite{Dwyer-Wilkerson}, this kind of gluing data is used to construct certain loop spaces with a free group action. To obtain the equivariant version of such a construction, one needs to work with the geometric realization of the enriched bar construction applied on an equvariant homotopy diagram. After forgetting the group action, one can only understand these constructions up to homotopy and in some cases, these actions cannot be realized. For example the symmetric group on three letters can not act freely on a sphere \cite{Milnor} however it can act freely on a finite CW complex homotopy equivalent to a sphere \cite{Swan}. For this reason, it is natural to ask when we can do such constructions up to homeomorphism. In this paper, we consider this question when the image of the homotopy diagram is a single point.

Given a small category $\mathcal{C}$,  let $\mathcal{F} \mathcal{C}$  (see Section \ref{subsect:standardresolution}) be the standard resolution of $\mathcal{C}$ considering $\mathcal{C}$ as a simplicially enriched category. It is known that $\mathcal{F} \mathcal{C}$ is a cofibrant replacement of $\mathcal{C}$ in the category of simplicially enriched categories with object set equal to the object set of $\mathcal{C}$. We define two functors $B$ and $\widetilde{B}$
from $\mathbf{Cat}$, the category of small categories, to $\mathbf{Top}$ which sends the category $\mathcal{C}$ to $|N\mathcal{C}|$, the geometric realization of the nerve of the category $\mathcal{C}$, and to $|B_{\bullet }(\ast,\mathcal{F} \mathcal{C} ,\ast)|$, the geometric realization of the enriched bar construction of $\mathcal{F} \mathcal{C}$, respectively. Using a cubical version of this simplicial resolution, namely $W\mathcal{C}$, (see Section \ref{subsect:Wconstruction}) as a middle step we prove that there is a homeomorphism from $\widetilde{B}\mathcal{C}=|B_{\bullet }(\ast,\mathcal{F} \mathcal{C} ,\ast)|$ to  $B\mathcal{C}=|N\mathcal{C}|$. One can ask whether the above homeomorphism is natural in $\mathcal{C}$, in other words, is there a natural homeomorphism between $\widetilde{B}$ and $B$. The following result partially answers this question.

\begin{thm}\label{thm:mainthm1}%[Theorem \ref{thm:doublebar-classifyinfspace}]
 There exists a natural homeomorphism between $\widetilde{B}$ and $B$ when they are restricted to the wide subcategory of $\mathbf{Cat}$ which contains all the inclusions.
\end{thm}

Since a group acts on a category by isomorphisms, the above result is especially important when one wants to use these constructions in the equvariant case. However to understand the maps and the homotopies of maps between equivariant constructions, we need a result that handles other morphisms. The category  $\mathbf{Pos}$, which is the full subcategory of $\mathbf{Cat}$ whose objects are posets considered as categories, is a good subcategory for this in the sense that the classifying space of any category is homeomorphic to the classifying space of a poset.

\begin{thm}\label{thm:mainthm2}%[Theorem \ref{thm:natural-Pos}]
There exists a twisted natural homeomorphism (Definition \ref{def:twistednattrans}) between $\widetilde{B}$ and $B$ when they are restricted to $\mathbf{Pos}$.
\end{thm}

Since the twisting of the twisted natural transformation in the above theorem is trivial on the wide subcategory in Theorem \ref{thm:mainthm1}, one can also use Theorem \ref{thm:mainthm2} for the equivariant case.

To obtain the above results, we introduce the notion of twisted natural transformation in Section \ref{sect:twistednaturaltransformations}. The key ingredient in the proof of Theorem \ref{thm:mainthm2}  is Theorem \ref{thm:extendingtwistednaturaltrans}, a  particular case of which is used to prove Theorem \ref{thm:mainthm1}. We consider it as one of our main results in this paper.

The paper is organized as follows: In Section \ref{sect:twistednaturaltransformations}, we introduce a natural transformation with a twist. This section can be considered as a summary of ideas behind the proofs of our main results. Section \ref{sect:definitionofBtilde} reminds the definition of simplicial sets, cubical sets, their realizations, some binary operations on such sets, standard resolution of a category, W-construction of a category, and finally the definition of enriched bar construction. These are all here for completeness and the reader familiar with these notations can skip this section. Section \ref{sect:tntbetweenBs} gives the statements and the proofs of Theorem \ref{thm:mainthm1} and Theorem \ref{thm:mainthm2}. Section \ref{sect:triangulation} and \ref{sect:Proofsofsomeresults} are devoted to proofs of technical lemmas which are used to prove our main theorems.

\section{Twisted natural transformation}
\label{sect:twistednaturaltransformations}

A natural transformation between two functors is usually defined as an assignment that assigns a morphism in the common codomain of the two functor to a given object in their common domain. An equivalent way to define a natural transformation can be given by assigning a morphism in their codomain to every morphism in their domain. Exploiting this equivalent way of defining a natural transformation here we define the notion of a natural transformation with a twist. Let $\mathcal{C}, \mathcal{M}$ be  categories, $\mathcal{I}$ be a wide subcategory of $\mathcal{C}$, and $F,G \maps \mathcal{C} \to \mathcal{M}$ be functors.

\begin{defn}\label{def:twistednattrans}A pair $(\alpha, t)$ is called a natural transformation from $F$ to $G$ with a twist on $G$ away from $\mathcal{I}$ if to every morphism $f \maps c \to c'$ in $\mathcal{C}$, $\alpha$ assigns a morphism $\alpha(f)\maps F(c) \to G(c')$ in $\mathcal{M}$ and $t$ assigns a morphism $t(f) \maps G(c) \to G(c)$ in $\mathcal{M}$ such that
\begin{enumerate}
  \item[i)] $\alpha(hgf)=G(h)t(h)\alpha(g)F(f)$ for every ternary composition $hgf$ in $\mathcal{C}$,
  \item[ii)] $t(f)=\mathrm{id}_c$ for every morphism $f\maps c\to c'$ in $\mathcal{I}$, and
  \item[iii)] $G(f_1)t(g_1)=t(g_2)G(f_1)$ for every $f_1,f_2 \in \mathcal{I}$ and $g_1, \ g_2 \in \mathcal{C}$ with $f_2g_1=g_2f_1$.
%  \begin{center}
%\begin{tikzpicture}[node distance=2cm,auto]
%\node (C) {$G(c)$};
%\node (A) [right of=C] {};
%\node (E) [right of=A] {$G(c)$};
%\node (D) [below of=C] {$G(c')$};
%\node (B) [below of=E] {$G(c').$};
%\draw[->] (C) to node {$t(g_1)$} (E);
%\draw[->] (C) to node [swap] {$G(f_1)$} (D);
%\draw[->] (E) to node {$G(f_1)$} (B);
%\draw[->] (D) to node {$t(g_2)$} (B);
%\end{tikzpicture}
%\end{center}
\end{enumerate} In particular, when $\alpha(\mathrm{id}_x)$ is an isomorphism for every object $x$ in $\mathcal{C}$, we say that $(\alpha, t)$ is a natural isomorphism with a twist on $G$ away from $\mathcal{I}$. In this case, we call $F$ a right twisted version of $G$.
\end{defn}

Given a functor $i\maps \mathcal{D}\to \mathcal{C}$, let  $(i/c)$ be the over category, i.e, the category whose objects are pairs $(d, \sigma)$  where $d$ is an object in $\mathcal{D}$ and $\sigma\maps i(d)\to c$ is a morphism in $\mathcal{C}$ and morphisms from $(d_1, \sigma_1)$ to $(d_2, \sigma_2)$ are morphisms  $f:d_1\to d_2$ in $\mathcal{D}$ so that $\sigma _1=\sigma _2 \circ i(f)$. We denote the projection to the first component by $\pi _c:(i/c) \to \mathcal{D}$. A functor $i \maps \mathcal{D}\to \mathcal{C}$ is said to be dense if for every object $c$ in $\mathcal{C}$  the colimit of $i\pi_c$ is naturally isomorphic to $c$.  We define a stronger version of being dense as follows. Let  $\mathcal{I}$ be a wide subcategory of $\mathcal{C}$. Set $(i/c)_{\mathcal{I}}$ to be the subcategory of $(i/c)$ so that an object $(d, \sigma)$ of $(i/c)$ is in  $(i/c)_{\mathcal{I}}$ if and only if $\sigma $ is in $\mathcal{I}$ and a morphism $g$ in $(i/c)$ is in  $(i/c)_{\mathcal{I}}$ if and only if $i(g)$ is in $\mathcal{I}$. We denote the inclusion $(i/c)_{\mathcal{I}}$ in $(i/c)$ by $i_c$.

\begin{defn}  A functor $i\maps \mathcal{D}\to \mathcal{C}$ is called  $\mathcal{I}$-dense (for $F:\mathcal{C}\to \mathcal{M}$) if the following conditions are satisfied
\begin{enumerate}
\item[i)] The functor $i \maps \mathcal{D}\to \mathcal{C}$ is dense (for $F$).
\item[ii)] For every object $c$ in $\mathcal{C}$, there exists a retract $r _c\maps (i/c) \to (i/c)_{\mathcal{I}}$ with a natural transformation (in the usual sense) $s$ from $1_{(i/c)}$ to $i_{c}r_{c}$  which is identity on $(i/c)_{\mathcal{I}}$. We write $s(d, \sigma)=s_{\sigma}$ and $r_{c}(d,\sigma)=(dr_{c}(\sigma), r_c(\sigma))$ for short.
\item[iii)] For every morphism $f: x \to y$ in $\mathcal{C}$, every object $(d,\sigma)$ in $(i/x)$, the morphism
 $r_y(f_*(s_{\sigma }))$ in $(i/y)$ is the identity  where $f_{\ast} \maps (i/x) \to (i/y)$ is the functor induced by $f$.
\end{enumerate}
\end{defn}
Here the condition iii) means that  $dr_y(f\sigma)=dr_yf(r_x(\sigma))$, $r_y(fr_x(\sigma))=r_y(f\sigma)$ and $r_{y}(f_{\ast}(s\sigma))=1_{dr_y(f\sigma)}$. Since $s$ is a natural transformation for any $\sigma\maps i(d) \to x$ in $i/x$ and$f \maps x \to y$ and $g \maps y \to z$ in $\mathcal{C}$, the following diagram
\begin{eqnarray*}\begin{CD} (d,gf\sigma)@>s_{gf\sigma}>> (dr_z(gf\sigma), r_z(gf\sigma))\\
@Vg_{\ast}(s_{f\sigma})VV@VVr_z(g_{\ast}(s_{f\sigma}))=\mathrm{id}V\\
(dr_y(f\sigma),gr_y(f\sigma))@>s_{gr_y(f\sigma)}>> (dr_z(g r_y(f\sigma)))
\end{CD}\end{eqnarray*}
commutes and hence we have $i(s_{gf\sigma})=i(s_{gr_y(f\sigma)})i(s_{f\sigma})$ in $\mathcal{C}$.

Note that if $i\maps \mathcal{D}\to \mathcal{C}$ is  full and dense for $F,G \maps \mathcal{C}\to \mathcal{M}$ then  a natural transformation from $F\circ i \maps \mathcal{D} \to \mathcal{M}$ to $G\circ i\maps \mathcal{D} \to \mathcal{M}$ induces a natural transformation from $F$ to $G$  using universal properties of colimits. Now we discuss a twisted version of this fact.

\begin{thm}\label{thm:extendingtwistednaturaltrans} Let $i:\mathcal{D}\to \mathcal{C}$ be full and $\mathcal{I}$-dense for $F,G \maps \mathcal{C} \to \mathcal{M}$. Then a natural transformation $(\alpha, t)$ from $F|_{\mathcal{D}}$ to $G|_{\mathcal{D}}$ with a twist on $G$ away from $\mathcal{I}\cap i(\mathcal{D})$ induces a natural transformation $(\bar{\alpha}, \bar{t})$ from $F$ to $G$ with a twist on $G$ away from $\mathcal{I}$ such that $(\bar{\alpha}, \bar{t})$ restricted to $i(\mathcal{D})$ is equal to $(\alpha, t )$.
\end{thm}
\begin{proof}  For a morphism $f\maps x \to y$ in $\mathcal{C}$, we define a morphism $\bar{\alpha }(f)$ from $F(x)$ to $G(y)$ by defining a natural transformation $A(f)$ from $Fi\pi_x$ to $\Delta _{G(y)}$ as follows:
$$A(f)_{(d,\sigma )}=G(r_y(f\sigma))\circ \alpha(i s_{f\sigma })\maps Fi\pi_x((d,\sigma ))\to G(y)$$
where $(d,\sigma )$ is an object in $(i/x)$. Here $\Delta_{G(y)}$ is a category with one object $G(y)$ and one morphism.  Similarly we define $\bar{t}(f)$ from $G(x)$ to $G(x)$ by defining a natural transformation $T$ from $Gi\pi_xi_x$ to $\Delta _{G(x)}$ as follows:
$$T(f)_{(d,\sigma )}=G(\sigma)\circ t(is_{f\sigma })\maps Gi\pi_x(i_x(d,\sigma ))\to G(x)$$
where $(d,\sigma )$ is an object in $(i/x)_{\mathcal{I}}$. For $f \in \mathcal{I}$, the map $is_{f\sigma}$ is in $I \cap i(\mathcal{D})$ and hence $\bar{t}(f)=\mathrm{id}$, i.e, $\bar{t}$ satisfies the condition ii).

Let $f\maps x\to y$, $g\maps y\to z$, and $h\maps z\to w$ be morphisms in $\mathcal{C}$. For any $(d,\sigma )$ object in $(i/x)$, we have
$s_{hgf\sigma}=s_{hr_z(gf\sigma)}\circ s_{gr_z(f\sigma)}\circ s_{f\sigma}$. We also have $G(r_z(gf\sigma))t(is_{hr_z(gf \sigma)})=\bar{t}(h)G(r_z(gf\sigma))$ and $\bar{\alpha}(g)F(r_y(f \sigma))=G(r_z(gr_y(f\sigma)))\alpha(is_{gr_y(f\sigma)})$ by construction. Therefore we have
\begin{eqnarray}
A(hgf)_{(d,\sigma )} &=& G(r_w(hgf\sigma)) \alpha (is_{hgf\sigma })  \nonumber \\
   &=& G(r_w(hgf\sigma)) G(i s_{hr_z(gf\sigma)}) t(i s_{hr_z(gf\sigma)}) \alpha (is_{gr_u(f\sigma ) })F(is_{f\sigma })
   \nonumber \\
   &=&  G(h)G(r_z(gf\sigma))t(is_{hr_z(gf\sigma)}) \alpha (is_{gr_y(f\sigma ) })F(is_{f\sigma })
   \nonumber \\
   &=&  G(h) \bar{t}(h)G(r_z(gf\sigma))\alpha (is_{gr_y(f\sigma ) })F(is_{f\sigma })\nonumber \\
   &=& G(h) \bar{t}(h)G(r_z(gr_y(f\sigma)))\alpha (is_{gr_y(f\sigma ) })F(is_{f\sigma }) \nonumber \\
   &=&  G(h) \bar{t}(h) \bar{\alpha }(g)F(r_y(f\sigma))F(is_{f\sigma })\nonumber \\
   &=& G(h) \bar{t}(h) \bar{\alpha }(g)F(f)F(\sigma ). \nonumber
\end{eqnarray}
 It yields the equality $$\bar{\alpha }(hgf)=G(h) \bar{t}(h) \bar{\alpha }(g)F(f),$$ that is, $(\bar{\alpha}, \bar{t})$ satisfies the first condition.

To show that condition iii) also holds for $\bar{t}$, let $\sigma \in (i/x)_{\mathcal{I}}$, $f_1\maps x \to y$, $f_2\maps z \to w$ be morphisms in $\mathcal{I}$ and $g_1 \maps x \to z$ and $g_2 \maps y \to w$ be morphisms in $\mathcal{C}$ with $g_2 \circ f_1=f_2 \circ g_1$. Then, we have $s_{g_2r_y(f\sigma)}s_{f_1\sigma}=s_{f_2r_z(g_1\sigma)}s_{g_1\sigma}$ and $s_{f_1\sigma}=s_{g_1 \sigma}=\mathrm{id}$. Therefore we have
\begin{eqnarray*} G(f_1)T(g_1)_{(d,\sigma)}&=& G(f_1)G(\sigma)t(is_{g_1\sigma})\\
&=& G(r_y(f_1\sigma))G(is_{f_1\sigma})t(is_{g_1\sigma})\\
&=& G(r_y(f_1\sigma))t(is_{g_2}r_y(f_1\sigma))G(is_{f_1\sigma})\\
&=& \bar{t}(g_2)G(r_y(f_1\sigma))G(is_{f_1\sigma})\\
&=& \bar{t}(g_2)G(f_1)G(\sigma).
\end{eqnarray*}

This proves that $(\bar{\alpha}, \bar{t})$ is a natural transformation from $F$ to $G$ with a twist on $G$ away from $\mathcal{I}$. It remains to show that  $\bar{\alpha}_{i(\mathcal{D})}=\alpha$ and $\bar{t}_{i(\mathcal{D})}=t$. Let $f \maps x \to y$ in $\mathcal{D}$. Since $i(\mathcal{D})$ is full $(d,\sigma) \in (i/i(x))$ if and only if $\sigma \in i(\mathcal{D})$. Therefore the equality $i(f)\sigma=r_y(i(f)\sigma)i(s_{i(f)\sigma})$ induces
\begin{eqnarray*} G(r_y(i(f)\sigma)\alpha(is_{i(f)\sigma})&=&\alpha(r_y(i(f)\sigma)i(s_{i(f)\sigma}))= \alpha(i(f)\sigma)
=\alpha(i(f))F(\sigma)
\end{eqnarray*} and hence $\bar{\alpha(i(f))}=\alpha(i(f))$.  On the other hand, the equality $\bar{t}(i(f))=t(i(f))$ directly follows from the fact that $(d,\sigma) \in (i/i(x))_{\mathcal{I}}$ if and only if $\sigma \in \mathcal{I} \cap i(\mathcal{D})$ and the last condition for $(\alpha,t)$ to be a natural transformation with a twist on $G$ away from $\mathcal{I} \cap i(\mathcal{D})$.
\end{proof}

\section{Definition of $\widetilde{B}$}
\label{sect:definitionofBtilde}

\subsection{Simplicial and cubical sets}
\label{subsect:simplicialandcubicalsets}
 The aim of this section is to introduce the necessary background on cubical sets. Although the theory of simplicial sets is more familiar, we also include the basic definitions of simplicial set theory for completeness. We refer reader to \cite{Goerss-Jardine} for more details on simplicial theory.

Let $[n]$ denote the category with objects $\{0,1,2,\dots ,n\}$
and exactly one morphism from $i$ to $j$ when $i<j$.
Let $\Delta$ be the finite ordinal number category, i.e.,
the full subcategory of the category $\mathbf{Cat}$ of small categories with objects
$[n]$ for $n$ a non-negative integer. A simplicial object in a small category $\mathcal{C}$ is a functor from
$\Delta ^{\op}$ to $\mathcal{C}$. Simplicial objects in $\mathcal{C}$
form a category denoted by $\mathbf{s}\mathcal{C}$ with morphisms given by natural
transformations. For $K$ a simplicial object in a category $\mathcal{C}$, we write $K_{n}$
instead of $K([n])$ and
call it the $n$-simplices of $K$.
For $i$ in $\{0,1,2,\dots n\}$, we define a functor $d^{i}$ from $[n-1]$ to $[n]$ by
$$d^{i}(j)=\begin{cases}
 j &\text{ if } j<i\\
j+1 &\text{ if }  j\geq i
\end{cases}$$
and define a functor $s^{i}$ from $[n+1]$ to $[n]$ by
$$s^{i}(j)=\begin{cases}
 j &\text{ if } j\leq i\\
j-1 &\text{ if } j> i.
\end{cases}
$$
It is know that all morphisms  of $\Delta $ can be written
as compositions of the above morphisms. For a simplicial object $K$,
the map
$d_{i}\maps K_{n}\to K_{n-1}$
induced by $d^{i}$ is called an $i$-th face map
of $K$, and the map $s_{i}\maps K_{n}\to K_{n+1}$ induced by $s^{i}$ is
called an $i$-th degeneracy map.

Objects of  $\mathbf{sSet}$ are called simplicial sets.
There is a standard embedding of $\Delta $ in $\mathbf{sSet}$  which takes $[n]$ to the standard $n$-simplex
$\Delta ^{n}=: \Delta(-, [n])\maps\Delta ^{\op}\to \mathbf{Set}$. Let $\Delta_{\leq n}$ be the full subcategory of the finite ordinal number category on the objects $\{[0], [1], \cdots, [n]\}$. The functor induced from the inclusion $\Delta_{\leq n} \hookrightarrow \Delta$ is called the $n$-truncation functor $\tau_n: \mathbf{Set}^{\Delta_{\leq n}^{\op}} \rightarrow \mathbf{sSet}$. This functor has a left adjoint and composing two gives the $n$-th skeleton functor $sk_n: \mathbf{sSet} \to \mathbf{sSet}$. Therefore, given a simplicial set $K$, $sk_nK$ is the subsimplicial set of $K$ generated
by $K_{0},\dots ,K_{n}$ under the degeneracies.

The category of simplicial sets is a symmetric monoidal category with the simplicial tensor product defined for simplicial sets $K$ and $L$ by
$$
(K\otimes L)_n~:=~K_n\times L_n~.
$$

The geometric realization of a simplicial set is a functor $|-|\maps \mathbf{sSet} \to \mathbf{Top}$ which is defined on objects by
$$|K|=~\colim_{\Delta^n \to K}|\Delta^n|=\Big( \coprod_{i \geq0}  K_n \times |\Delta^n| \Big) \Big/ \sim$$ where the equivalence relation has $(\phi^{\ast}x,t)=(x,\phi_{\ast}t)$ for every $\phi: [n] \to [m]$ in $\Delta$ and $(x,t) \in K_m\times |\Delta^{n}|$. Here the induced map $\phi_{\ast}\maps|\Delta^n|\to|\Delta^m|$ is defined by $\phi_{\ast}(t_0, \dots, t_n) = (s_0,\dots , s_m)$, where $s_i= \sum_{j \in \phi^{-1}(i)}t_j$. The geometric realization of a simplicial set is a $CW$-complex. Moreover, if we consider the geometric realization as a functor to the category $\mathbf{CGHaus}$ of compactly generated Hausdorff spaces then it preserves products.

For $i$ in $\{1,2,\dots n\}$ and $\varepsilon $ in $\{0,1\}$, define a function $d^{(i,\varepsilon )}$ from $I^{n-1}$ to $I^{n}$ by
$$d^{(i,\varepsilon )}(t_1,\dots, t_{n-1})=(t_1, \dots, t_{i-1}, \varepsilon, t_i, \dots, t_{n-1})$$
and define a function $s^{i}\maps I^{n}\to I^{n-1}$ by
$$s^{i}(t_1,\dots,t_n)=(t_1,\dots,t_{i-1},t_{i+1},\dots,t_n)
$$
and define a function $c^{i}$ from $I^{n+1}$ to $I^{n}$ by
$$c^{i}(t_1,\dots,t_{n+1})=(t_1,\dots,t_{i-1}, \mathrm{max}\{t_i,t_{i+1}\},t_{i+2},\dots,t_{n+1}).
$$
Let $\Box  _{c}$ denote the Box category with connections, whose objects are the $n$-cubes $I^n$, $n \geq 0$ and the morphisms of $\Box  _{c}$ are
generated by the functions $d^{(i, \varepsilon)}$, $s^{i}$ and
 $c^{i}$ where $n\geq 0$ and $1\leq i \leq n$.
Let $\mathcal{C}$ be a category, a functor $K\maps \Box  _{c}^{\op}\to \mathcal{C}$ is called a cubical object
of $\mathcal{C}$. We write $K_{n}$ instead of $K(I^{n})$ and
call it the $n$-cubes of $K$. The two maps
$d_{(i,\varepsilon)}\maps K_{n-1}\to K_{n}$
induced by $d^{(i,\varepsilon)}$ is called an $i$-th face map
of $K$, the map $s_{i}\maps K_{n}\to K_{n-1}$ induced by $s^{i}$ is
called an $i$-th degeneracy map, and the map $c_{i}\maps K_{n+1}\to K_{n}$ induced by $c^{i}$ is
called an $i$-th connection map.
The collection of all cubical objects of a category $\mathcal{C}$ form a category
which we denote by $\mathbf{c}\mathcal{C}$, in other words, $\mathbf{c}\mathcal{C}$
denotes the functor category from $\Box  _{c}^{\op}$ to $\mathcal{C}$.
Objects of $\mathbf{cSet}$ are called cubical sets.
There is a standard embedding of $\Box  _{c}$ in $\mathbf{cSet}$ called the standard $n$-cube
$\mathcal{I}^{n}$ which is the cubical set $ \Box  _{c}(-, I^{n})\maps\Box _{c}^{\op}\to \mathbf{Set}$.
For a cubical set $K$, the subcubical set of $K$ generated
by $K_{0},\dots ,K_{n}$ under the degeneracies is called the $n$-th
cubical skeleton of $K$, written $\skc _{n}K$.  As in the simplicial sets, one could also see the $n$-skeleton functor $\skc _{n}\maps\mathbf{cSet}\to\mathbf{cSet}$ as a composition of the $n$-truncation functor $\tau_{n}$ on cubical sets with its left adjoint.

The geometric realization of a cubical set is also defined similarly. It is a functor $|-|\maps \mathbf{sSet} \to \mathbf{Top}$ which is defined on objects by
$$|K|=~\colim_{\mathcal{I}^{n} \to K} I^n=\Big( \coprod_{i \geq0}  K_n \times I^n \Big) \Big/ \sim$$ where the equivalence relation has $(\phi^{\ast}x,t)=(x,\phi_{\ast}t)$ for every $\phi: [n] \to [m]$ in $\Delta$ and $(x,t) \in K_m\times I^{n}$.

The geometric realization does not preserves the levelwise product of cubical sets. For this reason, we work with the cubical product defined for cubical sets $K$ and $L$ by
$$
K\otimes L~:=~\colim _{\mathcal{I}^{j}\to K,~\mathcal{I}^{k}\to L}~I^{j+k}~.
$$  The product $I^n \otimes I^m=I^{n+m}$ defines a monoidal structure on the box category and the above product is indeed equal to the one induced from this monoidal product by Day convolution. In \cite{Jardine}, Jardine shows that the analogous monoidal product for cubical sets without connections agrees with the product structure on topological spaces. The same proof applies here. Indeed, since $|-|$ preserves colimits, it follows that
\begin{eqnarray*} |K \otimes L| & \cong & \colim_{\mathcal{I}^{n} \to K \times L} I^n \cong \colim_{\mathcal{I}^{j} \to K,~\mathcal{I}^{k} \to L}
I^{j+k} \\ & \cong & \colim_{\mathcal{I}^{j} \to K,~\mathcal{I}^{k} \to L} I^j \times I^k \cong |K| \times |L|.
\end{eqnarray*}
Here we consider the geometric realization as a functor to the category $\mathbf{CGHaus}$. Unfortunately, the monoidal product $\otimes$ on cubical sets is not symmetric. However since the monoidal category $(\mathbf{cSet}, \otimes)$ is biclosed, the category $\mathbf{cSet}$ is still enriched over itself.

There are adjoint functors
$$\mathbf{cSet}\rightleftharpoons^{T}_{S} \mathbf{sSet}$$ relating the cubical sets and the simplicial sets. The singular functor $S$ associates to every simplicial set $X$, a singular cubical set $S(X)$ with $n$-cells $S(X)_n=  \mathbf{sSet}((\Delta^1)^n, X)$. The triangulation functor $T$ is defined by $TK= \colim _{I^{n}\to K}(\Delta^1)^n$. These functors induce an equivalence of homotopy categories.

\subsection{Standard resolution of a category}
\label{subsect:standardresolution}
The standard resolution construction of \cite{Cordier} and \cite{Dwyer-Kan} (see also \cite{Steimle}) is a way of assigning a simplicial category to a small category. More precisely, the free category on a small category $\mathcal{C}$ is the free category $\mathcal{F} \mathcal{C}$, where
$\obj(\mathcal{F} \mathcal{C} )=\obj(\mathcal{C} )$,
and a morphism in $\mathcal{F} \mathcal{C} $ from $c$
to $d$ is a word of composable non-identity morphisms
$(f_n)\dots (f_1)$ such that the domain of $f_1$ is $c$ and the
range of $f_n$ is $d$. Here the composition is given by the concatenation
of words and the identity morphism on $c$ is the empty word $()$. For simplicity, we write $\mathcal{F}  f:=(f)\in\mathcal{F} \mathcal{C} (c,d)$
for a morphism $f\in\mathcal{C} (c,d)$.

Now, consider the functor $\mathcal{F} \maps \mathbf{Cat} \to \mathbf{Cat}$ which sends every category to its free category and which is defined on morphisms by sending a functor $F\colon\mathcal{C} \to\mathcal{D}$ to the induced one $\mathcal{F}  F\colon \mathcal{F} \mathcal{C} \to\mathcal{F} \mathcal{D} $ given by  $\mathcal{F}  f = \mathcal{F}  (F(f))$. It has a co-monad structure with the co-unit functor $U_\mathcal{C} \colon \mathcal{F} \mathcal{C} \to\mathcal{C} $
which is the identity on objects and sends $\mathcal{F}  f$ to $f$.
The co-multiplication functor
$\Delta_\mathcal{C} \colon \mathcal{F} \mathcal{C} \to \mathcal{F} ^2\mathcal{C} $
is defined by $\Delta_\mathcal{C}(\mathcal{F}  f)= \mathcal{F} ^2 f$. Clearly these functors satisfy the comonad identities.
Therefore one can construct a simplicial category (with discrete object set) $\mathcal{F} _\bullet \mathcal{C} $ by letting
$\mathcal{F} _n\mathcal{C} =\mathcal{F} ^{n+1}\mathcal{C} $.
Here the face and the degeneracy maps are identity on objects and on
morphisms they are given by
\begin{align*}
d_i=\mathcal{F} ^i(U_{\mathcal{F} ^{n+1-i}\mathcal{C} })& \colon
\mathcal{F} ^{n+1}\mathcal{C} (c,d)\to \mathcal{F} ^n\mathcal{C} (c,d)\\
s_i=\mathcal{F} ^i(\Delta_{\mathcal{F} ^{n-i}\mathcal{C} })& \colon
\mathcal{F} ^{n+1}\mathcal{C} (c,d)\to\mathcal{F} ^{n+2}\mathcal{C} (c,d).
\end{align*}

In other words, for $0\leq i \leq n-1$, the face map $d_i$ removes the parenthesis that are contained in exactly $i$ others while the degeneracy map $s_i$ doubles these parenthesis. When $i=n$, the map $d_n$ composes the innermost parenthesis and the map $s_n$ inserts the parenthesis to each map in the innermost parenthesis.

One can also consider $\mathcal{F} _{\bullet }\mathcal{C} $ as a simplicially enriched
category by taking objects as $\mathrm{obj}(\mathcal{C})$ and hom-sets as simplicial sets with $\mathcal{F}_{n}(c,d)=\mathcal{F}^{n+1}(c,d)$. Here what we mean by a simplicially enriched category is a category whose hom-sets are simplicial sets and composition maps are simplicial maps.  We refer reader to \cite{Kelly} for more information about enriched categories. One can consider the $0$-simplices as the maps and the $1$-simplices as the homotopies between maps, and the higher simplices as the higher homotopies between homotopies.
 A standard example of a simplicially enriched category is the category $\mathbf{Top}$ of compactly generated topological spaces which is enriched by defining $\mathbf{Top}(X,Y)$ by $$\mathbf{Top}(X,Y)_n= \mathbf{Top}(X \times |\Delta^n|, Y)$$
where $|\Delta^n|$ is the standard $n$-simplex in $\mathbb{R}^{n+1}$. Note that $\mathbf{Top}(X,Y)_{0}$ is the set of continuous maps from $X$ to $Y$.

A map between simplicially enriched categories consists of maps of $n$-simplices and called simplicial functor. More precisely, a simplicial functor $F\maps \mathcal{C} \to \mathcal{D}$ is given by a function $F: \mathrm{obj}(\mathcal{C}) \to \mathrm{obj}(\mathcal{D})$ together with simplicial maps $\mathcal{C}(c,c') \to \mathcal{D}(F(c),F(c'))$ which satisfy the functoriality diagrams. We denote the category of simplicially enriched categories and simplicial functors by $\mathbf{S}$-$\mathbf{Cat}.$

Now, we consider the standard resolution of the category $[n]$. From now on, we denote the map from $i$ to $i+1$ by $f_i$. It is well-known that the morphisms of the free simplicial resolution of $[n]$ is given as follows (See \cite{Cordier-Porter})
\begin{eqnarray*} \mathcal{F}_{\bullet}[n](i,j)= \left\{%
\begin{array}{ll}
    \Delta[1]^{j-i-1}, & \hbox{$i \leq j$;} \\
    \emptyset, & \hbox{otherwise.} \\
\end{array}%
\right.
\end{eqnarray*}
Here, we take $\Delta[1]^{-1}=\Delta[0]$.

\begin{ex}\label{ex:(0,3)}When $n=3$, $i=0$ and $j=3$, we have
\begin{eqnarray*} \mathcal{F}_0[3](0,3)&=& \{(f_2)(f_1)(f_0), (f_2)(f_1f_0), (f_2f_1)(f_0), (f_2f_1f_0)\}\\
\mathcal{F}_1[3](0,3)&=& \{((f_2))((f_1))((f_0)), ((f_2))((f_1f_0)), ((f_2f_1))((f_0)),((f_2f_1f_0)),
((f_2)(f_1)(f_0)),\\
&&((f_2)(f_1f_0)), ((f_2f_1)(f_0)), ((f_2)(f_1))((f_0)), ((f_2))((f_1)(f_0))\} \end{eqnarray*} where the first $4$ elements of $\mathcal{F}_2[3](0,3)$ are degenerate and others are non-degenerate. Moreover, $\mathcal{F}_{\bullet}[3](0,3)$ has exactly $2$ non-degenerate $2$-simplices $\sigma=(((f_2)(f_1))((f_0)))$ and $\tau= (((f_2))((f_1)(f_0)))$. The boundaries of these simplices are given as follows:
\begin{eqnarray*} d_0\sigma= ((f_2)(f_1))((f_0)) , \ d_1\sigma= ((f_2)(f_1)(f_0)) , \ d_2\sigma= ((f_2f_1)(f_0)),\\
d_0\tau=((f_2))((f_1)(f_0)) , \ d_1\tau=((f_2)(f_1)(f_0)) , \ d_2\tau= ((f_2)(f_1f_0)). \
\end{eqnarray*}

Since all the higher simplices are degenerate, $\mathcal{F}_{\bullet}[3](0,3)$ is the following simplicial set
{\tiny
\begin{center}
\begin{tikzpicture}[node distance=3cm,auto]
\node (C) {$(f_2)(f_1f_0)$};
\node (E) [right of=C] {};
\node (D) [below of=C] {$(f_2f_1f_0)$};
\node (A) [right of=E] {$(f_2)(f_1)(f_0)$};
\node (B) [below of=A] {$(f_2f_1)(f_0)$};
%\node (A) at (0,0) {$(f_2)(f_1)(f_0)$};
%\node (B) at (2,0) {$(f_2)(f_1f_0)$};
%\node (C) at (0,1) {$(f_2f_1)(f_0)$};
%\node (D) at (2,1) {$(f_2f_1f_0)$};
\draw[<-] (A) to node {$((f_2))((f_1)(f_0))$} (B);
\draw[<-] (A) to node [swap] {$((f_2)(f_1))((f_0))$} (C);
\draw[<-] (A) to node {} (D);
\draw[<-] (C) to node [swap] {$((f_2)(f_1f_0))$} (D);
\draw[<-] (B) to node  {$((f_2f_1)(f_0))$} (D);
\node () at (3,-1.5) {${\scriptsize((f_2)(f_1)(f_0))}$};
\node () at (2,-0.8) {{\Huge$\sigma$}}; % {\small(((f_2)(f_1))((f_0)))}
\node () at (4.2,-2.2) {{\Huge$\tau$}}; %{\small(((f_2))((f_1)(f_0)))}
\end{tikzpicture}
\end{center}}
\end{ex}

In general, all the simplices of $\mathcal{F}_{\bullet}[n](i,j)$ of dimension greater than $(j-i-1)$ are degenerate and it has $(j-i-1)!$ non-degenerate $(j-i-1)$-simplices faces of which give all the non-degenerate simplices of smaller dimensions. These $(j-i-1)!$ many $(j-i-1)$-simplices obtained by fully parenthesizing $f_{i},\cdots,f_{j-1}$ in the following way. The innermost parenthesis contains exactly one morphism $f_k$ of $[n]$, that is, $d_0^{j-i-1}=\underbrace{d_0 \cdots d_0}_{j-i-1}$ of the $(j-i-1)$-simplices are of the form $(f_{j-1})\cdots(f_i)$. Only one of the the second innermost parenthesis is of the form $((f_k)(f_{k+1}))$ and the others are of the form $((f_k))$. This means that $d_0^{j-i-2}$ of these simplices are of the form $$((f_i))\cdots((f_{k-1}))((f_k)(f_{k+1}))((f_{k+2}))\cdots((f_j)).$$ In general, $d_{0}^{j-i-s+1}$ of such a $(j-i-1)$-simplex is the composition of $j-i-s+2$ many morphisms in $\mathcal{F}_{s-2}[n](a,b)$ and $d_{0}^{j-i-s}$ is obtained by putting two consecutive such morphisms in the same parenthesis and inserting parenthesis around the other individuals.

Let $\sigma^n_{i,j} \maps  [j-i] \to [n]$ be defined by $\sigma^n_{i,j}(k)=i+k$. There is a one-to-one correspondence between the set of non-degenerate $(j-i-1)$-simplices of $\mathcal{F}_{\bullet}[n](i,j)$ and the set $L^n_{i,j}$ of chains $l=(\sigma_0,\cdots, \sigma_{j-i-1})$ where $\sigma_{j-i-1}=\sigma^n_{i,j}$ and $\sigma_{k}=d_{l_{k}}\sigma_{k+1}$ for some $ 1 \leq l_k \leq k+1$ ($l_0=1$). Here, the chain $l$ corresponds to the non-degenerate $(j-i-1)$-simplex $\tau_l=(d^0\tau_l)$ (since $l_0=1$) where $d_{0}^{j-i-s}(\tau_l)$ is obtained from $d_{0}^{j-i-s+1}(\tau_l)$ by taking the $l_{j-i-s}$-th and $(l_{j-i-s}+1)$-th morphisms in common parenthesis and inserting parenthesis around others. In the previous example $\sigma$ corresponds to the chain $l_{\sigma}=(d_1d_2\sigma^3_{0,3}, d_2\sigma^3_{0,3},\sigma^3_{0,3})$ and $\tau$ corresponds to the chain $l_{\tau}=(d_1d_1\sigma^3_{0,3}, d_1\sigma^3_{0,3},\sigma^3_{0,3})$.

A simpler way of representing $(j-i-1)$-simplex $\sigma_l$ of $\mathcal{F}_{\bullet}[n](i,j)$ is to look at the images of $\sigma_{k}$'s in $l$. Since  $\sigma_{k}=d_{l_{k}}\sigma_{k+1}$ means that the image of $\sigma_k$ is the set $\mathrm{Im}(\sigma_{k+1}) \backslash \{\sigma_{k+1}(l_k)\}$, we can consider $l$ as a tuple $(\tilde{l}_1,\cdots,\tilde{l}_{j-i-1})$ where $\mathrm{Im}(\sigma_{k+1})\backslash \mathrm{Im}(\sigma_{k})=\{\tilde{l}_{k+1}\}$. Therefore, we can represents $\sigma_l$ by this tuple. Let $\widetilde{L}^n_{i,j}$ be the set of $(j-i-1)$-tuples $(\tilde{l}_1,\cdots,\tilde{l}_{j-i-1})$ where $i+1 \leq \tilde{l}_k \leq j-1$ and $\tilde{l}_k \neq \tilde{l}_s$ if $k\neq s$.
Then there is a one-to-one correspondence between the set of non-degenerate $(j-i-1)$-simplices of $\mathcal{F}_{\bullet}[n](i,j)$ and the set $\widetilde{L}^n_{i,j}$. In the previous example $\sigma$ corresponds to the tuple $(1,2)$ while $\tau$ is corresponding to the tuple $(2,1)$.

%In the following sections, we interchangeably consider the geometric realization of $\mathcal{F}_{\bullet}[n](i,j)$ as a $(j-i-1)$ dimensional simplicial complex by using these correspondences.

\subsection{$W$-construction}
\label{subsect:Wconstruction}
In this section we introduce the cubical version of the standard resolution of a category which is due to Boardman and Vogt \cite{Boardman-Vogt}.

A cubically enriched category is a category whose hom-sets are cubical sets and composition maps are cubical maps. %{\color{red} \st{We refer reader to [9] for more information about enriched categories}.}
 One can consider the $0$-cubes of the hom-sets as the maps and the $1$-simplices as the homotopies between maps, and the higher cubes as the higher homotopies between homotopies. One standard example of a cubically enriched category is the category $\mathbf{Top}$ with
$$\mathbf{Top}_{\mathbf{cSet}}(X,Y)_n= \mathbf{Top}(X \times I^n, Y).$$ A map between cubically enriched categories consists of maps of $n$-cubes and called a cubical functor. More precisely a cubical functor $F \maps \mathcal{C} \to \mathcal{D}$ is given by a function $F \maps \mathrm{obj}(\mathcal{C}) \to \mathrm{obj}(\mathcal{D})$ together with cubical maps $\mathcal{C}(c,c') \to \mathcal{D}(F(c),F(d))$ which satisfy the functoriality diagrams. We denote the category of cubically enriched categories and cubical functors by $\mathbf{C}$-$\mathbf{Cat}$.

Given a small category $\mathcal{C}$, there is a specific way of assigning a cofibrant cubically enriched category $W \mathcal{C}$ called $W$-construction. Although it is first introduced by Boardman and Vogt topologically, here we give the more categorical version due to Blanc and Markl \cite{Blanc-Markl}. The $W$-construction $W\mathcal{C}$ on $\mathcal{C}$ is a cubically enriched category with $\mathrm{obj}(W \mathcal{C})=\mathrm{obj}(\mathcal{C})$ and for every $a,b$ objects of $\mathcal{C}$, the cubical set $W\mathcal{C}(a,b)$ constructed as follows:

$$ W\mathcal{C}(a,b)= \Big( \coprod_{\substack{
      \sigma\maps [n+1]\rightarrow \mathcal{C}, \\
      \sigma(0)=a, \ \sigma(n+1)=b}} \mathcal{I}_{\sigma}^n  \Big) \Big / \sim$$
where the equivalence relation is given by $$\mathcal{I}^{n-1}_{\sigma d^{n+1-i}} \sim d^{i,0}(\mathcal{I}^n_{\sigma}) \ \mathrm{and} \ \mathcal{I}^{n+1}_{\sigma s^{n+1-i}}\sim \begin{cases}
 s_i(\mathcal{I}^n_{\sigma}) &\text{ if } i=1, n+1\\
 c_i(\mathcal{I}^n_{\sigma})&\text{ otherwise.}
\end{cases}$$ The cubical composition
$$W\mathcal{C}(c,b)\otimes W\mathcal{C}(a,c)~\to~
W\mathcal{C}(a,b)$$
%$$
%W\mathcal{C}(a_{n+1-i},a_{n+1})\otimes W\mathcal{C}(a_{0},a_{n+1-i})~\to~
%W\mathcal{C}(a_{0},a_{n+1})=W\mathcal{C}(a,b)
%$$
sends the $(n-1)$-cube
$\mathcal{I}^{i-1}_{\beta}\otimes
       \mathcal{I}^{n-i}_{\sigma}$
to $d^{(i, 1)}\mathcal{I}^n_{\beta \cdot \sigma}$
where $(\beta \cdot \sigma)(j\to j+1)= \sigma(j \to j+1)$ if $j \leq n-i$ and $(\beta \cdot \sigma)(j \to j+1)= \beta((j+i-n-1)\to (j+i-n))$, otherwise. For example, for $\mathcal{C}=[n]$, the cubical set  $W[n](i,j)$ has only one non-degenerate $(j-i-1)$-cube corresponding to $\sigma^n_{i,j}$ and other cubes are glued to it as defined above and hence we have
\begin{eqnarray*} W[n](i,j)= \left\{%
\begin{array}{ll}
    \mathcal{I}^{j-i-1}, & \hbox{$i \leq j$;} \\
    \emptyset, & \hbox{otherwise.} \\
\end{array}%
\right.
\end{eqnarray*} Here, the cubical composition $W[n](j,k)\otimes W[n](i,j) \to W[n](i,k)$ is given by the cubical map $d_{k-j,1}\maps\mathcal{I}^{k-j-1} \times \mathcal{I}^{j-i-1} \to \mathcal{I}^{k-i-1}$ for $i\leq k\leq j$.

Moreover for every functor $f\maps \mathcal{C} \rightarrow \mathcal{D}$, the W-construction induces a cubical functor $Wf\maps W \mathcal{C} \to W \mathcal{D}$ which is given on objects by $f$ and on morphisms by sending an $n$-cube $\mathcal{I}^n_{\sigma}$ to $\mathcal{I}^n_{f \circ \sigma}$. This makes $W$ a functor from the category of small categories to the category of  cubically enriched categories.

For every small category $\mathcal{C}$, the simplicially enriched categories $TW\mathcal{C}$ and $\mathcal{F} _{\bullet }\mathcal{C}$ are naturally equivalent (see Lemma 3.6 in \cite{Baues-Blanc} and Remark 2.21 in \cite{Blanc-Markl}). Using this equivalence, one can consider $|W[n](i,j)|=I^{j-i-1}$ as a $(j-i-1)$-dimensional pure simplicial complex whose facets are $|\Delta_l^{j-i-1}|=[\sigma_0,\cdots, \sigma_{j-i-1}]$ where $l=(\sigma_0,\cdots,\sigma_{j-i-1}) \in L^n_{i,j}$. The cubical composition induces a map $c_{i,j,k} \maps |W[n](j,k)|\times |W[n](i,j)| \to |W[n](i,k)|$ which sends $(t'_1,\dots,t'_{k-j-1};t_1,\dots,t_{j-i-1})$ to $(t'_1,\dots,t'_{k-j-1},1,t_1,\dots,t_{j-i-1})$. Since the geometric realization of simplicial sets preserves the product, we can consider $|W[n](j,k)|\times |W[n](i,j)|$ as a pure simplicial complex whose facets are $(k-j-2)$-simplices  $[\sigma_{m_1}\times \tau_{m'_1},\dots,\sigma_{m_{k-i-1}} \times \tau_{m'_{k-i-1}}]$ where $m_s+m'_s=s-1$, $l=(\sigma_0,\cdots,\sigma_{k-j-1})$ in $L^n_{j,k}$ and $l'=[\tau_0,\cdots,\tau_{j-i-1}] \in L^n_{i,j}$. In this case, the map $c_{i,j,k}$ becomes a simplicial map which sends the vertex $\sigma_{m_s}\times \tau_{m'_s}$ to the vertex $\sigma_{m_s}\cdot \tau_{m'_s}$. Here, $c_{i,j,k}$ sends the simplex $[\sigma_{m_1}\times \tau_{m'_1},\dots,\sigma_{m_{k-i-1}} \times \tau_{m'_{k-i-1}}]$ to the simplex $[\beta_1,\dots, \beta_{k-j-1}]$ of $\Delta^{k-j-1}_{(\beta_0, \dots,\beta_{k-j-1})}$ where $\beta_s=\sigma_{m_s}\cdot \tau_{m'_s}$ for $1 \leq s \leq k-j-1$.

Under the homeomorphism $|W[n](i,j)|\cong I^{j-i-1}$, the vertex $\sigma_{i} \in l$ is sent to the vertex $(\epsilon_1,\cdots,\epsilon_{j-i-1})\in \mathbb{R}^{j-i-1}$ where $\epsilon_s=1$ if $(s+i)\in \mathrm{Im}(\sigma_i)$ and $\epsilon_{s}=0$, otherwise. For $x=(x_1, \dots, x_{k-j-1})$ in $I^{k-j-1}$ and $y=(y_1,\dots,y_{j-i-1})$ in $I^{j-i-1}$, let $x=\underset{m}{\sum}t_m \sigma_m \in |\Delta^{k-j-1}_{l}|$ and $y=\underset{m}{\sum}t'_m \tau_m \in |\Delta^{j-i-1}_{l'}|$ and  define  $z_s=x_s$ if $s \leq k-j-1$, $z_{k-j}=1$ and $z_s=y_{s-k+j+1}$, otherwise. Then $c_{i,j,k}(x,y)=\underset{m}{\sum}u_m \beta_m$ in $|\Delta^{k-i-1}_{l''}|$ with $\widetilde{l''}=(k-j+i,i_2+i,\cdots,i_{k-j-1}+i)$ and $u_s=z_{i_{s+1}}-z_{i_s}$ where $i_1=k-j$ and $1=z_{i_1} \geq z_{i_2} \geq \cdots \geq z_{i_{k-j-1}}$.

\subsection{Enriched bar construction}
\label{subsect:equivariantdoublebarconstruction}

We define the enriched bar construction as follows.

\begin{defn}Let $\mathcal{C}$ be a simplicially (respectively cubically) enriched category and
$F\maps \mathcal{C} \to \mathbf{Top}$ be a simplicial (respectively cubical) functor. The
enriched double bar construction
$B_{\bullet }(\ast, \mathcal{C} ,F)$ is a simplicial topological space with $n$-th simplices given as follows:
  \begin{equation*}B_{n }(\ast,\mathcal{C} ,F)=
    \coprod_{\begin{array}{c}
      \alpha_0, \alpha_1,\dots \alpha_n \\
       \text{ objects of }\mathcal{C}
    \end{array}}
   |\mathcal{C}(\alpha_{n-1},\alpha_n)\times
\dots \times \mathcal{C}(\alpha_{0},\alpha_1)|
    \times F(\alpha_0)
  \end{equation*}
where the face and the degeneracy maps are defined by using composition in
 $\mathcal{C}$, the evaluation maps $\mathcal{C}(\alpha_{0},\alpha_1) \times F(\alpha _0)\to F(\alpha_1)$
and the insertion of $\mathcal{C}(\alpha ,\alpha)$.
\end{defn}

The notion of a two-sided double bar construction for enriched categories was introduced by Shulman \cite{Shulman} and this construction can be considered as a one-sided version of Shulman's double bar construction. For this reason, we keep the Shulman's notation.

Now we are ready to define $\widetilde{B}$ as follows:
$$\widetilde{B}\mathcal{C}:=|B(\ast , \mathcal{F} \mathcal{C},\ast )|$$
equivalently we can define it as follows:
$$\widetilde{B}\mathcal{C}:=|B(\ast , W \mathcal{C},\ast )|$$

\section{Twisted natural transformation from $\widetilde{B}$ to $B$}
\label{sect:tntbetweenBs}

In this section, we will use the results of Section \ref{sect:twistednaturaltransformations} to construct a twisted natural transformation from $\widetilde{B}$ to $B$. For this we first study  $\widetilde{B}[n]$  and construct a homeomorphism from $\widetilde{B}[n]$ to $B[n]$. Then we give the twisting information on $B[n]$ for a given morphism $\sigma\maps [n]\to [m]$. We obtain our main results by working with appropriate categories. We postpone the proofs of all the lemmas given in this section to Section \ref{sect:Proofsofsomeresults} for the sake of legibility.

Recall that the morphisms of the $W$-construction of $[n]$ is given as follows
\begin{eqnarray*} W[n](i,j)= \left\{%
\begin{array}{ll}
    \mathcal{I}^{j-i-1}, & \hbox{$i \leq j$;} \\
    \emptyset, & \hbox{otherwise.} \\
\end{array}%
\right.
\end{eqnarray*} Define $B^n_0=\{0,1,\cdots,n\}$ and for $i>0$ $$B_i^n= \coprod_{r=(r_0, \cdots, r_i) \in R^n_i}|W[n](r_{i-1}, r_i)|\times |W[n](r_{i-2},r_{i-1})|\times \cdots \times
|W[n](r_{0}, r_1)|$$ where $R^n_i= \{r =(r_0, \cdots, r_i)| \ 0 \leq r_0< r_1<\cdots<r_i\leq n\}$. Then the geometric realization of the double bar construction $B(\ast, W[n], \ast)$ is given by
$$\widetilde{B}[n]= \Big( \coprod_{i \geq0}  B^n_i \times |\Delta^i| \Big) \Big/ \sim$$ where $(d_jx,t) \sim (x,d^jt)$ for every $(x,t) \in B^n_i\times |\Delta^{i-1}|$. From now on, we write $|C^n(r)|:=|W[n](r_{i-1}, r_i)|\times \cdots \times
|W[n](r_{0}, r_1)| $ for short. %We also denote an element $(t^{(i)}_1, \cdots, t^{(i)}_{r_i-r_{i-1}-1}; \cdots ; t^{(1)}_1, \cdots, t^{(1)}_{r_1-r_0-1}) $ of $|C^n(r)|$ by $t(r)$ and an element $(t_0, \cdots, t_i) \in |\Delta^i|$ by $t$.

To define a natural homeomorphism from $\widetilde{B}[n]$ to $B[n]$, we first define a map $\widetilde{f}[n](r) \maps |C_n(r)|\times |\Delta^i| \to |\Delta^n|$ as follows. For any $A \subseteq \{0,1,\cdots,n\}$, let $b_{A}=(b_0,b_1,\dots,b_n)$ where
$b_j=0$  if $j$ is not in $A$ and $b_j=1/|A|$ otherwise. For every morphism $\alpha \maps \{1, \cdots, i\} \to \{\pm 1\}$, let $|\Delta^i_{\alpha}|\subset |\Delta^i|$ be the $i$-simplex $$|\Delta^i_{\alpha}|=[b_{\alpha_0,\alpha_0}, \cdots, b_{\alpha_k,\alpha_k+k}, \cdots,b_{\alpha_i,\alpha_i+i}]$$ where $\alpha_{k}=|\{s \geq k+1| \ \alpha(s)=-1 \}|$. Here, we take $\alpha_i=0$. Let $U_i$ denote the set of all functions  $\alpha \maps \{1, \cdots, i\} \to \{\pm 1\}$. In Section \ref{sect:triangulation}, we show that the set $\{|\Delta^i_{\alpha}| \ | \ \alpha \in U_i\}$ gives a triangulation of $|\Delta^i|$.

 We also write $\{U_i\}=\underset{\alpha \in U_i}{\bigcup} \{(\alpha_k, \alpha_{k}+k)| \ 0\leq k \leq i \}$ and $\{L^n_{j,k}\}=\underset{l \in L^n_{j,k}}{\bigcup} \{l\}.$  Here $\{l\}=\{\sigma_0,\cdots \sigma_{j-i-1}\}$ where $l=(\sigma_0,\cdots,\sigma_{j-i-1})$. Then we define
 $$f[n](r)\maps \{L^n_{r_{i-1},r_i}\}\times \cdots \times \{L^n_{r_{0},r_1}\}\times \{U_i\}  \to \mathcal{P}(\{0,1,\dots,n\})$$
 by sending the tuple $(\sigma^i,\cdots,\sigma^1,(\alpha_k,\alpha_k+k))$ to the intersection of the union of images of $\sigma^j$'s and the set $\{r_{\alpha_{k}},r_{\alpha_{k}}+1,\cdots, r_{\alpha_k+k}-1, r_{\alpha_k+k}\}$. For each $(l^i,\cdots,l^1, \alpha)$ in $L^n_{r_{i-1},r_i}\times \cdots \times L^n_{r_{0},r_1}\times U_i$, this induces a map $$\widetilde{f}^{n,r}_{(l^i,\cdots,l^1, \alpha)}\maps |\Delta^{r_i-r_{i-1}-1}_{l^i}|\times \cdots \times |\Delta^{r_1-r_{0}-1}_{l^1}| \times |\Delta^{i}_{\alpha}| \to |\Delta^n|$$  which sends the point
$$\Big( \sum_{m} t^i_{m}\sigma_m^i, \cdots, \sum_{m}t^1_{m}\sigma^1_m, \sum_{0\leq k \leq i}\tilde{t}_{k}b_{\alpha_k,\alpha_k+k}\Big)$$ to the point $$\sum t^i_{m_i}\cdots t^1_{m_1}\tilde{t}_{k}b_{f[n](r)(\sigma_{m_i}^i,\cdots,\sigma_{m_1}^1,  (\alpha_k,\alpha_k+k))}$$ where the sum  is taken over all $(\sigma_{m_i}^i,\cdots,\sigma_{m_1}^1, (\alpha_s,\alpha_s+s)) \in (l^i,\cdots,l^1, \{\alpha\})$.

\begin{ex} Let $\sigma$ and $\tau$ in $W[3](0,3)$ are defined as in the Example \ref{ex:(0,3)}. Let $l_{\sigma}=(\sigma_0, \sigma_1, \sigma_2)$ and $l_{\tau}=\tau_0,\tau_1,\tau_2$ be the corresponding chains in $L^3_{0,3}$. Then for $\alpha,\beta:\{1\} \rightarrow \{\pm 1\}$ with $\alpha(1)=1$ and $\beta(1)=-1$, we have
\begin{eqnarray*}
  \widetilde{f}^{3,(0,3)}_{(l_{\sigma}, \alpha)}\bigg(\sum_{i=0}^2t_i\sigma_i,\sum_{i=0}^1\widetilde{t_i}b_{0,i}\bigg) &=& \bigg(\widetilde{t_0}+\bigg(\frac{t_0}{2}+\frac{t_1}{3}+\frac{t_2}{4}\bigg)\widetilde{t_1}, \bigg(\frac{t_1}{3}+\frac{t_2}{4}\bigg)\widetilde{t_1},\frac{t_2}{4}\widetilde{t_1},\bigg(\frac{t_0}{2}+\frac{t_1}{3}+\frac{t_2}{4}\bigg)\widetilde{t_1}\bigg)  \\
  \widetilde{f}^{3,(0,3)}_{(l_{\tau}, \alpha)}\bigg(\sum_{i=0}^2t_i\tau_i,\sum_{i=0}^1\widetilde{t_i}b_{0,i}\bigg) &=& \bigg(\widetilde{t_0}+\bigg(\frac{t_0}{2}+\frac{t_1}{3}+\frac{t_2}{4}\bigg)\widetilde{t_1}, \frac{t_2}{4}\widetilde{t_1}, \bigg(\frac{t_1}{3}+\frac{t_2}{4}\bigg)\widetilde{t_1},\bigg(\frac{t_0}{2}+\frac{t_1}{3}+\frac{t_2}{4}\bigg)\widetilde{t_1}\bigg)
  \\
  \widetilde{f}^{3,(0,3)}_{(l_{\sigma}, \beta)}\bigg(\sum_{i=0}^2t_i\sigma_i,\sum_{i=0}^1\widetilde{t_i}b_{1-i,1}\bigg) &=&  \bigg(\bigg(\frac{t_0}{2}+\frac{t_1}{3}+\frac{t_2}{4}\bigg)\widetilde{t_1}, \bigg(\frac{t_1}{3}+\frac{t_2}{4}\bigg)\widetilde{t_1},\frac{t_2}{4}\widetilde{t_1},\widetilde{t_0}+\bigg(\frac{t_0}{2}+\frac{t_1}{3}+\frac{t_2}{4}\bigg)\widetilde{t_1}\bigg)
  \\
  \widetilde{f}^{3,(0,3)}_{(l_{\tau}, \beta)}\bigg(\sum_{i=0}^2t_i\tau_i,\sum_{i=0}^1\widetilde{t_i}b_{1-i,1}\bigg) &=& \bigg(\bigg(\frac{t_0}{2}+\frac{t_1}{3}+\frac{t_2}{4}\bigg)\widetilde{t_1},\frac{t_2}{4}\widetilde{t_1}, \bigg(\frac{t_1}{3}+\frac{t_2}{4}\bigg)\widetilde{t_1},\widetilde{t_0}+\bigg(\frac{t_0}{2}+\frac{t_1}{3}+\frac{t_2}{4}\bigg)\widetilde{t_1}\bigg).
  \end{eqnarray*}
\end{ex}

Since $\widetilde{f}^{n,r}_{(l^i,\cdots,l^1, \alpha)}$'s are obtained by using the vertices of the corresponding simplices, we can glue these maps to obtain a continuous map $\widetilde{f}[n](r) \maps |C_n(r)| \times |\Delta^i| \to |\Delta^n|$.

\begin{lem}\label{lem:equivalencerespecting} For $(x,t) \in |C_n(r)|\times |\Delta^{i-1}|$ we have $$\widetilde{f}[n](d_jr)(d_jx,t)=\widetilde{f}[n](r)(x,d^jt)$$ where $d_jr=[r_0, \dots, r_{j-1},r_{j+1},\dots,r_{i}]$.
\end{lem}
The above lemma says that the maps $\widetilde{f}[n](r)$ respect the equivalence relations discussed above and hence they form a map $\widetilde{f}_n \maps \widetilde{B}[n] \to \Delta^n$. This maps also satisfies the following result.
\begin{lem}\label{lem:naturallityofn} For all $n>0$, the following diagram commutes
\begin{eqnarray*}\begin{CD}
\widetilde{B}[n-1]@>d^j_{\ast}>> \widetilde{B}[n] @. \\
@V\widetilde{f}_{n-1}VV@V\widetilde{f}_{n}VV  \\
B[n-1] @>d^j>> B[n],@. \\
\end{CD}
\end{eqnarray*}
that is, $d^j\widetilde{f}_{n-1}=\widetilde{f}_nd^j_{\ast}$.
\end{lem}
Let $\mathbf{Cat_i}$ denote the wide subcategory of $\mathbf{Cat}$ which has inclusions as morphisms.  Then the morphisms in $\mathbf{Cat_i} \cap \Delta$ is generated by the face maps $d^j\maps [n-1] \to [n]$. Hence the above result shows that $\widetilde{f}_n$ induces a natural transformation from the functor  $\widetilde{B} \circ i \maps \mathbf{Cat_i} \cap \Delta \to \mathbf{Top}$ to the functor $B \circ i \maps \mathbf{Cat_i} \cap \Delta \to \mathbf{Top}$.
Now we show that this natural transformation is in fact a natural homeomorphism.
\begin{lem}\label{lem:homeomorphism} For each $n$, the map $\widetilde{f}_n$ is a homeomorphism.
\end{lem}
Since the inclusion functor $i \maps \mathbf{Cat_i} \cap \Delta \to \mathbf{Cat_i}$ is full and dense, the following theorem immediately follows from Theorem \ref{thm:extendingtwistednaturaltrans}.
\begin{thm}\label{thm:doublebar-classifyinfspace}
There is a natural homeomorphism from $\widetilde{B}$ to $B$ restricted to $\mathbf{Cat}_i$.
\end{thm}
To generalize this result to category of posets, we take $\mathcal{I}$ to be the intersection of $\mathbf{Cat_i}$ and $\mathbf{Pos}$
 and $\mathcal{D}$ to be $\Delta$. We also introduce a twisting on $B[n]$ as follows. For each  morphism $f \maps [n]\to [m]$ in $\Delta $, we define $t(f )\maps \Delta ^n\to \Delta ^n$ as the affine extension of:
$$t(f )(b_A)=\sum _{i\in A}\frac{1}{|f(A)||A\cap f^{-1}(\{f(i)\})|}\, b_{\{i\}}$$
where $A$ is a subset of $\{0,1,2,\dots ,n\}$. Here $b_A=(b_0, \cdots, b_n)$ where $b_i=1/|A|$ if $i \in A$, $b_i=0$ otherwise.  Then we have the following key lemma
\begin{lem}\label{lem:twistswork}
Let $h$ and $t$ be defined as above, Then $(h, t)$ is a natural homeomorphism from $\widetilde{B}|_{\Delta}$ to $B|_{\Delta}$ with a twist on $B|_{\Delta}$ away from $\mathbf{Cat}_{i}\cap \Delta $
\end{lem}
Using this lemma we obtain the following result
\begin{thm}\label{thm:natural-Pos} There exists a natural homeomorphism $(\bar{\alpha}, \bar{t})$ from $\widetilde{B}$ to $B$ with a twist on $B$ restricted to $\mathbf{Pos}$ away from $\mathbf{Cat}_{i}\cap \mathbf{Pos}$.
\end{thm}

\section{A triangulation of $|\Delta^n|$}
\label{sect:triangulation}

The aim of this section is to introduce a specific triangulation of $|\Delta^n|$ as mentioned in the previous section, which is also crucial for the proofs given in the next section.

For $\alpha \maps \{1, \cdots, n\} \to \{\pm 1\}$, let $|\Delta^n_{\alpha}|$ be defined as in the previous section. Also define the sign of $\alpha$ by $\mathrm{sign}(\alpha)=(j_1, \cdots, j_m)$ where  $\{j_1,\cdots, j_m\}$ is the set of all integers satisfying $\alpha(j_k) \neq \alpha(j_k-1)$ with the order $1 \leq j_1 < \cdots < j_m\leq n$.

\begin{lem} \label{lem:delta-alpha} If $\alpha(1)=+1$ then the point $(x_0, \cdots, x_n) \in |\Delta^n|$ belongs to the $n$-simplex $|\Delta^n_{\alpha}|$ if and only if it satisfies the following inequalities
\begin{eqnarray*} x_{n}+\cdots + x_{\alpha_0+1} &\leq & x_0 + \cdots+ x_{\alpha_0}, \\
x_0+\cdots+x_{\alpha_{j_{2k+1}}} &\leq& x_n + \cdots + x_{(\alpha_{j_{2k+1}}+j_{2k+1})},\\
x_n + \cdots +x_{(\alpha_{j_{2k}}+j_{2k})} &\leq& x_0+\cdots + x_{\alpha_{j_{2k}}}
\end{eqnarray*}
for $1 \leq 2k,2k+1 \leq m$. On the other hand, if $\alpha(1)=-1$ then a point $(x_0,\cdots,x_n)$ in the standard $n$-simplex belongs to $|\Delta^n_{\alpha}|$ if and only if it satisfies the following inequalities
\begin{eqnarray*} x_{n}+\cdots + x_{\alpha_0} &\geq & x_0 + \cdots + x_{\alpha_0-1}, \\
x_0+\cdots+x_{\alpha_{j_{2k+1}}} &\geq& x_n + \cdots + x_{(\alpha_{j_{2k+1}}+j_{2k+1})},\\
x_n + \cdots +x_{(\alpha_{j_{2k}}+j_{2k})} &\geq& x_0+\cdots + x_{\alpha_{j_{2k}}}
\end{eqnarray*} for $1 \leq 2k,2k+1 \leq m$.
\end{lem}
\begin{proof} A point $x=(x_0, \cdots, x_n)$ belongs to $|\Delta^n_{\alpha}|$ if and only if $x= \sum_{k=0}^{n} a_k b_{\alpha_k,\alpha_k+k}$ for some $0 \leq a_k \leq 1$ with $\sum_{i=0}^na_i=1$. Since the cases $\alpha(1)=1$ and $\alpha(1)=-1$ are dual, we assume that $\alpha(1)=1$. In this case, we have
\begin{eqnarray*} x_{\alpha_0}&=&a_0+\frac{1}{2}\big (a_{1}+\cdots +a_{j_1-1}\big)\\
x_{\alpha_{j_{2k-1}}+j_{2k-1}} &=& \frac{1}{2} \big( a_{j_{2k-1}-1}+\cdots +a_{j_{2k}-1}\big)\\
x_{\alpha_{j_{2k}}}&=& \frac{1}{2}\big(a_{j_{2k}-1}+ \cdots + a_{j_{2k+2}-1}\big)\\
x_u&=& \left\{%
\begin{array}{ll}
\frac{1}{2}a_{\alpha_{j_{2k+2}}+j_{2k+2}-1-u}, & \hbox{$\alpha_{j_{2k+2}}<u<\alpha_{j_{2k}}$;} \\
    \frac{1}{2}a_{\alpha_{j_{2}}+j_{2}-1-u}, & \hbox{$\alpha_{j_2}<u<\alpha_0$;} \\
    \frac{1}{2}a_{u-\alpha_0}, & \hbox{$\alpha_0< u <\alpha_{j_1}+j_1$;} \\
            \frac{1}{2}a_{u-j_{2k}}, & \hbox{$\alpha_{j_{2k-1}+j_{2k-1}}< u <\alpha_{j_{2k+1}}+j_{2k+1}$.} \\
 \end{array}%
\right.\\
\end{eqnarray*} Then the result follows immediately.
\end{proof}

 Let $V=\{b_{i,j}| \ i \leq j \}$ and $$S=\bigcup_{\alpha \in U_n}\{ \mathrm{ \ set \ of \ all \ subsets \ of} \ \{b_{\alpha_0,\alpha_0}, \cdots, b_{\alpha_k,\alpha_k+k}, \cdots,b_{\alpha_n,\alpha_n+n}\}\}$$ where $U_n=\{\alpha \maps \{1, 2, \cdots, n\} \to \{\pm 1\}\}$. Then by the following lemma $(V,S)$ gives a triangulation of $|\Delta^n|$ with $2^n$ many $n$-simplices.

\begin{lem} $\bigcup_{\alpha \in U_n}|\Delta^n_{\alpha}|=|\Delta^n|.$
\end{lem}
\begin{proof} For each $x=(x_0, \dots, x_n) \in |\Delta^n|$, one can choose $\alpha$ for which $x \in |\Delta^n_{\alpha}|$ as follows. If $x_n \leq x_0$, let $\alpha(n)=+1$, if not let $\alpha(n)=-1$. Without loss of generality, let $x_n \leq x_0$. Then take $\alpha(j)=+1$ for each $n\geq j \geq i_1$ where $$x_n+ \cdots + x_{i_1} \leq x_0 < x_n+\cdots +x_{i_1-1}.$$
Then let $i_2$ be the largest integer for which
$$x_n+ \cdots + x_{i_1} \leq x_0+ \cdots +x_{i_2} < x_n+\cdots + x_{i_1-1}.$$
Then let $\alpha(j)=-1$ for each $ i_2-i_1-1 \leq j \leq i_1-1$. By continuing this way, we can determine $\alpha$ when $x_n \leq x_0$. This proves the lemma.
\end{proof}

Let $E_n$ be the subspace of the standard simplex $|\Delta^n|$ whose points $x=(x_0, \cdots, x_n)$ satisfy the relation $\mathrm{min}(x_i,x_j)\leq x_{i+1}$ for all $i<j$. Note that $E_n$ is indeed the topological space $A^{n}_{(0,1,\cdots,n)}$ defined in the following section. The subspace $E_i$ can be considered as the pure subcomplex of the barycentric subdivision of $|\Delta^n|$ who has $2^n$-many $n$-simplices of the barycentric subdivision of the standard simplex. More precisely, $E_n$ is the union of $n$-simplices with vertices
$\{b_{A_0},b_{A_1},\dots , b_{A_n}\}$ where $A_0$ has one element,
$A_n=\{0,1,2,\dots n\}$, $A_{j}=A_{j+1}- \min A_{j+1} $ or
$A_{j}=A_{j+1}- \max A_{j+1}$. The above triangulation of $|\Delta^n|$ also has $2^n$-many $n$-simplex and we can define a homeomorphism between $|\Delta^n|$ and $E^n$ as follows.
For each $\alpha$ and $1 \leq j \leq n$, let $A_j^{\alpha}$ be defined by the relations $A_n^{\alpha}=\{0, \cdots, n\}$ and
\begin{eqnarray*}A_j^{\alpha}= \left\{%
\begin{array}{ll}
    A_{j+1}^\alpha - \mathrm{min} \ A^{\alpha}_{j+1}, & \hbox{$\alpha(j+1)=-1$;} \\
A_{j+1}^\alpha - \mathrm{max} \ A^{\alpha}_{j+1}    , & \hbox{otherwise.} \\
\end{array}%
\right. \end{eqnarray*}
Then the map $f_n^{\alpha}$ defined by sending the vertex $b_{\alpha_j, \alpha_{j}+j}$ to the vertex  $b_{A^{\alpha}_j}$ gives a homeomorphism between $|\Delta^n_{\alpha}|$ and $ [b_{A^{\alpha}_0}, \dots, b_{A^{\alpha}_n}]$. Since the homeomorphisms $f_n^{\alpha}$ respect all the intersections, we can glue them together to obtain a homeomorphism $f_n \maps |\Delta^n| \to E_n$ as desired.

\section{Proofs of Lemmata in Section \ref{sect:tntbetweenBs}}
\label{sect:Proofsofsomeresults}

This section is devoted to proofs of the lemmas of Section \ref{sect:tntbetweenBs}. %We will consider $|C_n(r)| \times |\Delta|^i$ as a $(r_i-r_0)$-dimensional pure simplicial complex whose $(r_i-r_0)$-dimensional simplices are the facets of the subcomplexes of the form $|\Delta^{r_i-r_{i-1}}_{l^i}|\times \cdots \times |\Delta^{r_1-r_{0}}_{l^1}| \times |\Delta^{i}_{\alpha}|$ where $(l^i,\cdots,l^1, \alpha) \in L^n_{r_{i-1},r_i}\times \cdots \times L^n_{r_{0},r_1}\times U_i$. Let $l^j=(\sigma^j_0,\dots,\sigma^j_{r_j-r_{j-1}-1})$.
We first prove Lemma \ref{lem:equivalencerespecting}. For this, note that for $1\leq j \leq i-1$, we have $$d_j(x)=(x^i, \dots, x^{j+2}, c_{r_{j-1},r_{j},r_{j+1}}(x^{j+1},x^{j}), x^{j-1}, \dots,x^1)$$  where  $x=(x^i, \cdots,x^1) \in |\Delta^{r_i-r_{i-1}-1}_{l^i}|\times \cdots \times |\Delta^{r_1-r_{0}-1}_{l^1}|$. Here, $c_{r_{j-1},r_{j},r_{j+1}}(x_{j+1},x_{j})$ is in one of the $(k$-$j$-$2)$-simplices $[\beta_1,\dots, \beta_{r_{j+1}-r_{j-1}-2}]$ where $\beta_s=\sigma^{j+1}_{m_s}\cdot \sigma^j_{m'_s}$ for $1 \leq s \leq k-j-1$ and $m_s+m_s'=s-1$.
Moreover, for $t=\underset{0\leq k \leq i-1}{\sum}t_k b_{\alpha_k,\alpha_k+k}$ in $|\Delta^{i-1}_{\alpha}|$ we have $$d^jt=\underset{0\leq i \leq n}{\sum}t_i d^jb_{\alpha_i,\alpha_i+i}=\underset{0\leq k \leq i}{\sum} t_k' b_{\alpha'_k} \in |\Delta^i_{\alpha'}|$$ where
\begin{eqnarray*} \alpha_k'&=&\left\{
                                                                                        \begin{array}{ll}
                                                                                          \alpha_k+1, & \hbox{$k \leq i-1$;} \\
                                                                                          0, & \hbox{$k=i$.}
                                                                                        \end{array}
                                                                                      \right. ,  \ \mathrm{and} \
t_k'=\left\{
            \begin{array}{ll}
            t_k, & \hbox{$k\leq i-1$;} \\
            0, & \hbox{$k=i$.}
            \end{array}
            \right., \ \mathrm{for} \ j=0, \\
\alpha_k'&=&\left\{
                                                \begin{array}{ll}
                                                \alpha_k+1, & \hbox{$k \leq \nu_j$;} \\
                                                \alpha_{k-1}, & \hbox{$k>\nu_j$.}
                                                \end{array}
                                                \right. ,  \ \mathrm{and} \
t_k'=\left\{
       \begin{array}{ll}
         t_k, & \hbox{$k<\nu_j$;} \\
         0, & \hbox{$k=\nu_j$;} \\
         t_{k-1}, & \hbox{$k>\nu_j$.}
       \end{array}
     \right., \ \mathrm{for} \ 0<j\leq \alpha_0,\\
 \ \alpha_k'&=&\left\{
                                                \begin{array}{ll}
                                                \alpha_k, & \hbox{$k \leq \omega_j$;} \\
                                                \alpha_{k-1}, & \hbox{$k> \omega_j$.}
                                                \end{array}
                                                \right. ,  \ \mathrm{and} \
t_k'=\left\{
       \begin{array}{ll}
         t_k, & \hbox{$k<\omega_j$;} \\
         0, & \hbox{$k=\omega_j$;} \\
         t_{k-1}, & \hbox{$k>\omega_j$.}
       \end{array}
     \right. \ \mathrm{for} \ \alpha_0<j\leq i-1\\
  \end{eqnarray*} where $\nu_j$ is the smallest integer in $\{0,1, \cdots, i-1\}$ with $\alpha_{\nu_j}=j-1$ and $\omega_j \in \{0,1, \cdots, i-1\}$ is the smallest one with $\alpha_{\omega_j}+\omega_j=j$.

\begin{proof}[Proof of Lemma \ref{lem:equivalencerespecting}] The cases $j=0$ and $j=i$ are clear. For $0<j<i$, let $x=(x^i,\dots,x^1)$ be in $|\Delta^{r_i-r_{i-1}-1}_{l^i}|\times \cdots \times |\Delta^{r_1-r_{0}-1}_{l^1}|$ with $x^j= \underset{m}{\sum}t^j_m \sigma_m^j$ and $t= \underset{0\leq s \leq i-1}{\sum}\tilde{t}_{s}b_{\alpha_s,\alpha_{s}+s} \in |\Delta^{i-1}_{\alpha}|$. Then $d^jt= \underset{0\leq k \leq i}{\sum} t_k' b_{\alpha'_k} \in |\Delta^i_{\alpha'}|$ where $t'_k$ and $\alpha'$ are defined as above. Moreover $d_jx=(\tilde{x}^{i-1}, \dots, \tilde{x}^1)$ where $\tilde{x}^{k}=x^{k+1}$ if $k\geq j+1$, $\tilde{x}^j=\sum t_{m_{j+1}}^{j+1}t_{m_j}^j \big(\sigma_{m_{j+1}}^{j+1}\cdot \sigma_{m_j}^j\big)$ and $\tilde{x}^k=x^{k}$ for $k<j$.

For $(\alpha_s,\alpha_{s}+s) \in \{\alpha\}$, let $B_{r,\alpha_s}=\{r_{\alpha_s},r_{\alpha_s}+1,\dots,r_{\alpha_s+s}\}$. Then we have
\begin{eqnarray*}B_{r,\alpha's}&=&\left\{
                                    \begin{array}{ll}
                                      \{r_{\alpha_s+1}, r_{\alpha_s+1}+1,\dots, r_{\alpha_s+s+1}\}, & \hbox{$s \leq \nu_j, \ \mathrm{and} \ j \leq \alpha_0$;}\\
                                      \{r_{\alpha_{s-1}}, r_{\alpha_{s-1}}+1,\dots, r_{\alpha_{s-1}+s}\}, & \hbox{$s > \nu_j, \ \mathrm{and} \ j \leq \alpha_0$;}\\
\{r_{\alpha_s}, r_{\alpha_s}+1,\dots, r_{\alpha_s+s}\}, & \hbox{$s \leq \omega_j, \ \mathrm{and} \ j > \alpha_0$;}\\
\{r_{\alpha_{s-1}}, r_{\alpha_{s-1}}+1,\dots, r_{\alpha_{s-1}+s}\}, & \hbox{$s > \omega_j, \ \mathrm{and} \ j > \alpha_0$.}\\
                                    \end{array}
                                  \right.\\ \end{eqnarray*}
and
\begin{eqnarray*}
B_{d_jr,\alpha_s}&=&\{(d_jr)_{\alpha_s},\dots,(d_jr)_{\alpha_s+s}\}\\
&=&\left\{
                                    \begin{array}{ll}
                                      \{r_{\alpha_s}, r_{\alpha_s}+1,\dots, r_{\alpha_s+s}\}, & \hbox{$\alpha_{s}+s < j$;}\\
                                      \{r_{\alpha_{s}}, r_{\alpha_{s}}+1,\dots, r_{\alpha_{s}+s+1}\}, & \hbox{$\alpha_s < j < \alpha_s+s$;}\\
                                     \{r_{\alpha_{s}+1}, r_{\alpha_{s}+1}+1,\dots, r_{\alpha_{s}+s+1}\}, & \hbox{$\alpha_s \geq j $;}\\
\end{array}
                                  \right.\\
&=&\left\{
                                    \begin{array}{ll}
                                      \{r_{\alpha_s+1}, r_{\alpha_s+1}+1,\dots, r_{\alpha_s+s+1}\}, & \hbox{$s \leq \nu_j, \ \mathrm{and} \ j \leq \alpha_0$;}\\
                                      \{r_{\alpha_{s}}, r_{\alpha_{s}}+1,\dots, r_{\alpha_{s}+s+1}\}, & \hbox{$s > \nu_j, \ \mathrm{and} \ j \leq \alpha_0$;}\\
\{r_{\alpha_{s}}, r_{\alpha_{s}}+1,\dots, r_{\alpha_{s}+s}\}, & \hbox{$s \leq \omega_j, \ \mathrm{and} \ j > \alpha_0$;}\\
\{r_{\alpha_s}, r_{\alpha_s}+1,\dots, r_{\alpha_s+s+1}\}, & \hbox{$s > \omega_j, \ \mathrm{and} \ j > \alpha_0$.}\\
                                    \end{array}
\right.\\
\end{eqnarray*}since $s\leq \nu_j$ and $j \leq \alpha_0$ if and only if $\alpha_s \geq j$; $s>\nu_j$ and $j \leq \alpha_0$ if and only if $\alpha_s<j\leq \alpha_0< \alpha_s+s$; $s \leq \omega_j$ and $j> \alpha_0$ if and only if $\alpha_s+s<j$ ; and finally $s>\omega_j$ and $j> \alpha_0$ if and only if $\alpha_s+s>j<\alpha_0>\alpha_s$.

Since the image of $\sigma_{m_{j+1}}^{j+1}\cdot \sigma_{m_{j}}^{j}$ equals to the union of images of $\sigma_{m_{j+1}}^{j+1}$ and $\sigma_{m_{j}}^{j}$ and the coefficients $t'_{\nu_j}=t'_{\omega_j}=0$ and $t'_k=t_k$ if $k<\nu_j$ (or $k<\omega_j$) and $t'_k=t_{k-1}$ when $k>\nu_j$ (or $k> \omega_j$), the result easily follows from the above equalities.
\end{proof}

We now give a proof for Lemma \ref{lem:naturallityofn}.\\

\begin{proof}[Proof of Lemma \ref{lem:naturallityofn}] The induced map $d^j_{\ast}\maps \widetilde{B}[n-1] \rightarrow \widetilde{B}[n]$ is explicitly given as follows. Let $(x^i, \cdots, x^1,t) \in |\Delta^{r_i-r_{i-1}-1}_{l^i}|\times \cdots \times |\Delta^{r_1-r_{0}-1}_{l^1}|\times |\Delta^i_{\alpha}|$ be such that $$x^k=\underset{0\leq m \leq r_k-r_{k-1}-1}{\sum}t^k_m\sigma_m^k.$$ Then $d_{\ast}^j(x^i, \cdots, x^1,t)=(\overline{x}^i, \cdots, \overline{x}^1,t)$ where $$\overline{x}^k=\underset{0\leq m \leq r'_k-r'_{k-1}-1}{\sum}\overline{t}^k_m\overline{\sigma}_m^k$$ with
\begin{eqnarray*} r_k'=\left\{%
\begin{array}{ll}
    r_k, & \hbox{$j>r_{k}$,} \\
    r_k+1, & \hbox{otherwise,}
\end{array}%
\right.
 \quad \overline{t}^k_m=\left\{%
 \begin{array}{ll}
    0, & \hbox{$r_{k-1}<j \leq r_k, \ m=r_{k}-r_{k-1}$,} \\
 t^k_m, & \hbox{otherwise,}
\end{array}%
\right.
\end{eqnarray*}
and
\begin{eqnarray*} \overline{\sigma}^k_m(s)=\left\{%
\begin{array}{ll}
    \sigma_m^k(s), & \hbox{$j>r_k$;} \\
\sigma^k_m(s), & \hbox{$j=r_k, \ s\neq m+1, \ m\neq r_k-r_{k-1}$;} \\
\sigma^k_m(s), & \hbox{$r_{k-1}< j<r_k, \ s\neq m+1, \ m\neq r_k-r_{k-1}, \ \sigma^k_m(s)<j$;} \\
\sigma^k_m(s)+1, & \hbox{$r_{k-1}< j<r_k, \ s\neq m+1, \ m\neq r_k-r_{k-1}, \ \sigma^k_m(s)\geq j$;} \\
    r_{k}+1, & \hbox{$r_{k-1} < j \leq r_k, \ s=m+1$;}\\
    s+r_{k-1}, & \hbox{$r_{k-1}<j \leq r_k, \ m= r_k-r_{k-1}$;}\\
    \sigma_m^k(s)+1, & \hbox{$j\leq r_{k-1}$.}
\end{array}%
\right.
\end{eqnarray*}
Using this formula, one can easily check that the map $h_n$ respects the face maps in the sense that it makes the first diagram above commute.
\end{proof}

In order to show that $\widetilde{f}$ is a homeomorphism, we first show that the image of $\widetilde{f}[n](r)$ is indeed the subset $A_r^n$ of $B[n]$ consisting of all $a=(a_0,\cdots,a_n)$ satisfying the following conditions
\begin{itemize}
  \item[(1)] $a_m=0$ when $m> r_i$ or $m<r_0$,
  \item[(2)] $a_u \leq \mathrm{min}(a_{r_j},a_{r_{j+1}})$ when $r_j< u < r_{j+1}$,
  \item[(3)] $\mathrm{min}(a_{r_p}, a_{r_s}) \leq a_{r_{p+1}}$  when  $p<s$,
\end{itemize}
for each $r=(r_0, \cdots, r_i) \in R^n_i$ and $1 \leq i \leq n$.
%If $V$ is a set of vertices of a simplex $|\Delta^d|$, we write $|\Delta^d|=:|\Delta V|$. For a function $f \maps V_1 \times \cdots \times V_n \to W$ where %$V_i$ and $W$ are sets, let $$\widetilde{f}\maps|\Delta V_1| \times \cdots \times |\Delta V_m| \to |\Delta W|$$ be the simplicial map defined by %$$\widetilde{f}\big( \sum_{v \in V_1} t_{1v}v, \cdots, \sum_{v \in V_n}t_{nv}v\big)=\sum_{(v_1,\dots,v_n)\in V_1 \times \cdots \times V_n} t_{1v_1}\cdots %t_{nv_n}f(v_1,\cdots,v_n).$$

Recall that the barycentric subdivision of the standard $n$-simplex $|\Delta^n|$ is the pure simplicial complex with $n!$ many $n$-simplices with vertices
$\{b_{A_0},b_{A_1},\dots , b_{A_n}\}$ where $A_0\subset A_1\subset \cdots \subset A_n=\{0,1,2,\dots n\}$. Note that the facet $[b_{A_0},\cdots, b_{A_n}]$ of the barycentric subdivision corresponds to the subspace of all points $(x_0,\cdots,x_n)$ with $x_{k^A_0} \geq x_{k^A_1} \geq \cdots \geq x_{k^A_n}$ where $\{k^A_0\}=A_0$ and $\{k^A_j\}=A_j \backslash A_{j-1}$ otherwise. If we consider $|\Delta^{r_i-r_0}|$ as a subcomplex of $|\Delta^n|$ whose points are of the form $(0, \dots, x_{r_{0}}, \cdots,x_{r_i},0,\dots,0)$, then $A_r^n$ can be considered as the pure subcomplex of the barycentric subdivision of $|\Delta^{r_i-r_0}|$ which is the union of facets satisfying the conditions $(2)$ and $(3)$ above.

For $\alpha \in U_i$, let $A^{\alpha}=(A_0^{\alpha},\dots,A^i_{\alpha})$ be defined as in Section \ref{sect:triangulation} and for $l^j \in L^n_{r_{j-1},r_j}$, let $\widetilde{l^j}=(\tilde{l}^j_1,\dots, \tilde{l}^j_{r_j-r_{j-1}-1})$ in $\widetilde{L}^n_{r_{j-1},r_j}$ be defined as in Section \ref{subsect:standardresolution}. Then we have the following lemma.

\begin{lem} The image of $\widetilde{f}^{n,r}_{(l^i,\cdots,l^1, \alpha)}$ is the set of all points $(x_0,\cdots,x_n)$ of $|\Delta^n|$ satisfying the conditions (1), (2) and the following relations
\begin{itemize}
  \item[(2$'$)] $x_{\tilde{l}^j_1} \geq \cdots \geq x_{\tilde{l}^j_{r_j-r_{j-1}-1}}$ for $1 \leq j \leq i$,
  \item[(3$'$)] $x_{r_{k_0^{A^{\alpha}}}} \geq \cdots \geq x_{r_{k_i^{A^{\alpha}}}}$.
\end{itemize}
\end{lem}

\begin{proof} The lemma follows from the fact that  $$f(\sigma^i_{m_i},\cdots,\sigma^1_{m_1},(\alpha_k,\alpha_k+k))=\{r_{\alpha_k},r_{\alpha_k+1},\dots,r_{\alpha_k+k},\tilde{l}^{s+1}_1,\dots,\tilde{l}^{s+1}_{m_{s+1}},\dots,\tilde{l}^k_1,\dots,
\tilde{l}^k_{m_k}\}$$ where $l^j=(\sigma^j_0,\cdots,\sigma^j_{r_j-r_{j-1}})$.
\end{proof}

Note that the condition (3$'$) in the above lemma implies the condition (3) and the condition (2) and (2$'$) tells us how to order $x_u$'s satisfying condition (2) for fixed $k$.  We are now ready to prove \ref{lem:homeomorphism}. Let $A^{n,r}_{(l^i,\cdots,l^1, \alpha)}$ be the subset of $A^{n,r}$ satisfying 2$'$ and 3$'$.

\begin{proof}[Proof of Lemma \ref{lem:homeomorphism}] Let $(b_0, \cdots, b_n)$ be an point in $|\Delta^n|$. Let $r=(r_0, \cdots, r_{i})$ where $r_0$ is the smallest integer for which $b_u \neq 0$, $r_{i}$ is the largest integer for which $b_u\neq 0$, and $r_k$ is defined inductively as follows.
Define $r_k$ be the unique integer such that $b_{r_k}>b_{r_{k-1}}$ and $b_{r_{k-1}} \geq b_u$ for all $r_{k-1}<u<r_{k}$. This procedure terminates at $t$ when $b_u \leq b_{r_t}$ for all $u>r_t$. If $r_t=r_i$, we are done. If not, define $r_{i-k}$'s inductively as the unique integer such that $ b_{r_{i-k}} >b_{r_{i-k+1}}$  and $b_{r_{i-k+1}}\geq b_u$ for all $r_{i-k} < u < r_{i-k+1}$. This procedure terminates when $r_{i-k}=r_t$. We have $b_{r_0} < b_{r_1}< \cdots < b_{r_t}> b_{r_{t+1}}>\cdots> b_{r_i}$. By construction  $r=(r_0, r_1, \cdots, r_i)$ is the smallest subset of $\{1, \cdots, n\}$ for which $(b_0,\cdots, b_n) \in A^n_{r}$. Indeed if $(b_0,\cdots, b_n) \in A^n_{r'}$ where $r'=(r_0',\cdots,r'_{i'})$ and $r_k \notin \{r_0',\cdots,r'_{i'}\}$ then there is $p$ for which $r'_p < r_k< r'_{p+1}$ and hence $b_{r_k}\leq \mathrm{min}\{b_{r_p}, b_{r_{p+1}}\}$. However, we have $b_{r_k}> b_u$ for all $u<r_k$ if $k\leq t$ and $b_{r_k}> b_u$ for all $u>r_k$ if $k>t$ by construction. So $r_k \in \{r_0',\cdots,r'_{i'}\}$. Similarly one can show that if $(b_0, \cdots, b_n)$ is in the relative interior of $A^n_r$ then $r$ is uniquely determined.

 Now we define $l^1,\dots l^i,\alpha $ so that $(b_0,\dots,b_n)$ lies in the image of $\widetilde{f}^{n,r}_{(l^i,\cdots,l^1, \alpha)}$.  Let  $b_{r_{j_0}} \geq \dots \geq b_{r_{j_i}}$.  Define $\alpha \in U_i$ by the relations $\alpha_s=\mathrm{min}\{j_0,\cdots,j_s\}$. Similarly, we can choose $l^j$ by looking at the ordering of $a_u$'s for $r_{j-1} < u < r_{j}$. More precisely, if $a_{u_1} \geq \cdots \geq a_{u_{r_{j}-r_{j-1}-1}}$ for $ r_{j-1} < u_m < r_j$ then we define $l^j$ by $\tilde{l}^j_{m}=u_{m}$. Then $(b_0,\dots,b_n)$ lies in the image of $\widetilde{f}^{n,r}_{(l^i,\cdots,l^1, \alpha)}$. This means $\widetilde{f}_n$ is surjective. As above, one can show that if $(b_0,\cdots,b_n)$ is in the relative interior of $A^{n,r}_{(l^i,\cdots,l^1, \alpha)}$ then $l_1, l_2, \dots l_i, \alpha$ are uniquely defined. Note that a point in the image of $\widetilde{f}^{n,r}_{(l^i,\cdots,l^1, \alpha)}$ is an element of the relative interior of $A^{n,r}_{(l^i,\cdots,l^1, \alpha)}$ if and only if its preimage lies in the interior of $|\Delta^{r_i-r_{i-1}-1}_{l^i}|\times \cdots \times |\Delta^{r_1-r_{0}-1}_{l^1}| \times |\Delta^i_{\alpha}|$.  Since each $\widetilde{f}^{n,r}_{(l^i,\cdots,l^1, \alpha)}$ is injective restricted to the interior of its domain, it follows by induction that the map $\widetilde{f}_n$ is injective. Notice that the domain of $\widetilde{f}_n$ is compact and its codomain is Hausdorff. Hence $\widetilde{f}_n$ is open. This means $\widetilde{f}_n$ is a homeomorphism.

\end{proof}

\begin{proof}[Proof of Lemma \ref{lem:twistswork}]
It is straight forward to show that the conditions in Definition \ref{def:twistednattrans} are satisfied when we only consider the vertices of the barycentric subdivision. On rest of the domain of this function, it is defined by affine extensions hence the conditions in Definition \ref{def:twistednattrans} are satisfied on whole the domain.
\end{proof}

\end{document}